\documentclass[12pt,oneside]{amsart} 

\usepackage{hyperref}
\usepackage{fullpage}
\usepackage{graphpap}
\usepackage{graphicx}
\usepackage{epstopdf}
\usepackage{url}
\usepackage{amssymb}
\usepackage[left=1in,top=0.8in,right=0.8in,bottom=0.8in]{geometry}

\newcommand{\comments}[1]{}
\newcounter{parcount}[subsubsection]
\newtheorem{theorem}{Theorem}[section]

\newtheorem{lemma}[theorem]{Lemma}
\newtheorem{corollary}[theorem]{Corollary}
\DeclareGraphicsExtensions{.eps}

\newtheorem{utv*}{Proposition}
\newtheorem{hyp*}{Conjecture}

\newtheorem*{th*}{Theorem}

\numberwithin{equation}{section}

\newcommand{\av}[2]{\left\langle #1\right\rangle_{_{\scriptstyle #2}}}

\def\B{\mathcal{B}}

\numberwithin{equation}{section}

\theoremstyle{plain}

\title[Mutual estimates for $RH_p^d$ and $A_q^d$ constants for the doubling weights.]{Mutual estimates for the dyadic Reverse H\"{o}lder and Muckenhoupt constants for the dyadically doubling weights.}

\author{O. Beznosova}
\address{Department of Mathematics, Baylor University, One Bear Place \#97328, Waco, TX 76798-7328, USA.}
\address{Oleksandra\_Beznosova@baylor.edu}

\author{T.~Ode}
\address{Temitope\_Ode@baylor.edu}

\begin{document}

\begin{abstract}
Muckenhoupt and Reverse H\"{o}lder classes of weights play an important role in harmonic analysis, PDE's and quasiconformal mappings. In 1974 Coifman and Fefferman showed that a weight belongs to a Muckenhoupt class $A_p$ for some $1<p<\infty$ if and only if it belongs to a Reverse H\"{o}lder class $RH_q$ for some $1<q<\infty$. In 2009 Vasyunin found the exact dependence between $p$, $q$ and the corresponding characteristic of the weight using the Bellman function method. The result of Coifman and Fefferman works for the dyadic classes of weights under an additional assumption that the weights are dyadically doubling. We extend the Vasyunin's result to the dyadic Reverse H\"{o}lder and Muckenhoupt classes and obtain the dependence between $p$, $q$, the doubling constant and the corresponding characteristic of the weight. More precisely, given a dyadically doubling weight in $RH_p^d$ on a given dyadic interval $I$, we find an upper estimate on the average of the function $w^{q}$ over the interval $I$. From the bound on this average we can conclude, for example, that $w$ belongs to the corresponding $A_{s_1}^d$ class or that $w^p$ is in $A_{s_2}^d$ for some values of $s_i$. We obtain our results using the method of Bellman functions.
\end{abstract}

\maketitle

\section{Definitions and main results.}

We will be dealing with a family of dyadic intervals on the real line:
$$
D := \{[n 2^{-k}, (n+1)2^{-k}], \;\; n, k \in \mathbb{N} \}.
$$
For an interval $J$, let $D(J)$ stand for the family of all its dyadic subintervals, $D(J):= \{I \in D, I \subset J\}$ and let $D_n(J)$ stand for the family of all dyadic subintervals of $J$ of the length exactly $2^{-n}|J|$.
For a locally integrable function $f$, let $\av{f}{I}$ stand for the average of $f$ over the interval $I$, $\av{f}{I} := \frac{1}{|I|} \int_I f(x) dx$ where $|I|$ is the Lebesgue measure of $I$.

Let $w$ be a weight, i. e. $w$ is a locally integrable almost everywhere positive function. 
Since we will be dealing mostly with averages, we define the dyadic doubling constant of the weight $w$ to be
$$
Db^d (w) := \sup_{I\in D} \frac{\av{w}{I^\ast}}{\av{w}{I}},
$$
where $I^\ast$ is the dyadic ``parent'' of the interval $I$, i.e. the smallest dyadic interval that strictly contains the interval $I$. If the dyadic doubling constant of the weight $w$ is bounded by $Q$, we will say that $w \in Db^{d,Q}$. Note also that any weight is positive almost everywhere, therefore the dyadic doubling constant defined this way is always greater than $\frac{1}{2}$.

Our main assumption is that a weight $w$ belongs to the dyadic Reverse H\"older class of weights on the interval $J$ with the corresponding constant bounded by $\delta$:
$$
w \in RH_p^{\delta, d} (J) \;\;\;\;\;\text{if and only if}\;\;\;\;\;\; [w]_{RH_p^{\delta, d} (J)} := \sup_{I \in D(J)} \frac{\av{w^p}{I}^{1/p}}{\av{w}{I}} \leqslant \delta.
$$ 

We define the $A_{q}^{\delta, d} (J)$ to be the class of the dyadic Muckenhoupt weights on the interval $J$ with the corresponding constant bounded by $\delta$:
$$
w \in A_q^{\delta, d} (J) \;\;\;\;\;\text{if and only if}\;\;\;\;\;\; [w]_{A_q^{\delta, d} (J)} := \sup_{I \in D(J)} \av{w}{I} \av{w^{-\frac{1}{q-1}}}{I}^{q-1} \leqslant \delta.
$$
Given a dyadically doubling weight $w \in RH_p^{\delta,d} (J)$, our goal in this paper is to bound 
averages involved in definitions of $w \in A_q^d$ and $w^p \in A_q^d$:
$$
\av{w}{J} \av{w^{-\frac{1}{q-1}}}{J}^{q-1} \;\;\;\; \text{and} \;\;\;\; \av{w^p}{J} \av{w^{-\frac{p}{q-1}}}{J}^{q-1}.
$$
Note that quantities $\av{w}{J}$ and $\av{w^p}{J}$ are involved in the definition of $RH_p^{\delta,d} (J)$, therefore for our goals it is enough to bound from above $\av{w^s}{J}$ for $s<0$.

We are ready to define the Bellman function for our problem: for $p>1$, $s<0$ and $Q> \frac{1}{2}$, let
$$
\B (x_1, x_2; p, s, \delta, Q) :=  \sup_{w\in RH_{p}^{\delta, d} (J) ,\, Db^d (w) \leqslant Q}\{ \av{w^s}{J}: \;\; w \; is \; s.t.\;  \av{w}{J} = x_1, \av{w^{p}}{J} = x_2\}.
$$
When choice of the parameters $p$, $s$, $\delta$ and $Q$ is clear from the context, we will skip them and write $\B (x_1, x_2)$. Note also that by the rescaling argument, $\B$ does not depend on the interval $J$.

Then for the given $p$, $s$, $\delta$, and $Q$, $\B$ is defined on the domain
$$
\Omega_{\delta} := \{ \vec{x} = (x_1,x_2) \; : \; \exists w \in RH_p^{\delta,d} \; s.t.\; Db^d(w) \leqslant Q \; \text{and}\; x_1 = \av{w}{J} , x_2 = \av{w^p}{J}  \}.
$$

In order to state the main theorem we need to define functions $u^\pm_p (t)$. Let $u^\pm_p (t)$ be two solutions (positive and negative) of the equation 
\begin{equation}
\label{eqn:defUpm}
(1-p u)^\frac{1}{p} (1-u)^{-1} = t , \;\;\;\; 0 \leqslant t\leqslant 1.
\end{equation}
For $Q\geqslant 2$, we define $\varepsilon (p, \delta, Q)$ as follows. Let $H:= H(p, Q) = \frac{Q^p - 1}{Q-1} $ and $\varepsilon := \frac{H}{p} \left(\frac{{p} - 1}{H-1}\right)^{\frac{p - 1}{p}} \delta$.

Then we can define $$s^{\pm} (\varepsilon) := u^\pm \left(\frac{1}{\varepsilon}\right)$$ and $$r^{\pm} =  u^\pm \left(\frac{y^{1/p}}{\varepsilon x}\right).$$

Note that since $u^+(t)$ is a decreasing function and in our domain $\frac{1}{\varepsilon} \leqslant \frac{y^{1/p}}{\varepsilon x}$, we have that $r^+ \in [0,s^+]$. Similarly, since $u^-(t)$ is an increasing function, $r^- \in [s^-,0]$.


\begin{theorem}[Main Theorem] \label{thm:Main}
If $p>1$ $Q \geqslant 2$ and $\delta > 1$, let $s^- := s^- (\varepsilon)$ for $\varepsilon (p, \delta, Q)$ defined above. 
$$
\text{If $q \in \left(\frac{1}{s^-}, 0\right)$ then}\;\;\;\;\;\;
\B(x, y; p, q, \delta) \leqslant x^{q} \frac{1-qr^-}{1-qs^-} \left(\frac{1- s^-}{1-r^-}\right)^{{q}}=y^{\frac{q}{p}} \frac{1-qr^-}{1-qs^-} \left(\frac{1-p s^-}{1-p r^-}\right)^{\frac{q}{p}}.
$$
\end{theorem}

Proof of Theorem \ref{thm:Main} can be found in Section \ref{s:2}.

Note that the result from \cite{Vasyunin:08} is assuming that the Reverse H\"older inequality for the weight $w$ holds for any interval $I\subset J$, while our Theorem \ref{thm:Main} only uses dyadic subintervals $I \in D(J)$ and the doubling constant. Therefore our result is more general. Unfortunately, we lose the sharpness.


As a consequence of Theorem \ref{thm:Main}, one can easily obtain the following corollary.

\begin{corollary}[$RH_p$ vs $A_q$] \label{thm:MainRHpAq}
 Let $w$ be a Reverse H\"older dyadically doubling weight with $[w]_{RH_p^d} = \delta$ and $Q := \max \{Db^d(w),2\}$. 
 Let $\varepsilon (p, \delta, Q)$ be defined as above. Let $s^- = s^- (\varepsilon)$. Then
 
(i) for every $q > 1-s^-$ $w \in A_q^d$, and moreover
$$
[w]_{A_q^d} \leqslant \left(\frac{q-1}{q-1+s^-}\right)^{q-1};
$$
(ii) for every $q>1-ps^-$ $w^p \in A_q^d$, and moreover
$$
[w^p]_{A_q^d} \leqslant \left(\frac{q-1}{q-1+ps^-}\right)^{q-1}.
$$
Where $s^-(\varepsilon)$ is the negative solution of the equation $(1-p s^-)^\frac{1}{p} (1-s^-)^{-1} = \frac{1}{\varepsilon}$

\end{corollary}

A result similar to the second part of the above corollary was used in \cite{BMP:2013} for the sharp norms of $t$-Haar multiplier operators.


\section{Proof of Theorem \ref{thm:Main}} \label{s:2}

In this section we essentially follow the proof of Lemma 2 from \cite{Vasyunin:08}. Unfortunately, we cannot use full proof from Vasyunin's paper since it relies on the Lemma 4 from his paper, which fails in the dyadic case.

Let
$$
\B (x_1, x_2; p, s, \delta, Q) :=  \sup_{w\in RH_{p}^{\delta, d} (J) ,\, Db^d (w) \leqslant Q}\{ \av{w^s}{J}: \;\; w \; is \; s.t.\;  \av{w}{J} = x_1, \av{w^{p}}{J} = x_2\}
$$
and 
$$
B_{\max} = B_{\max}(x_1, x_2; p, q, \delta, Q) := x_1^{{q}} \frac{1-qr^-}{1-qs^-} \left(\frac{1- s^-}{1- r^-}\right)^{{q}}
$$
be defined on domains 
$$
\Omega_{\delta} = \{ \vec{x} = (x_1,x_2) \; : \; \exists w \in RH_p^{\delta,d} \; s.t.\; Db^d(w) \leqslant Q \; \text{and}\; x_1 = \av{w}{J} , x_2 = \av{w^p}{J}  \}
$$
 and 
 $$
 \Omega_\varepsilon := \{ \vec{x} = (x_1,x_2) \; : \; x_i > 0, \;x_1^p \leqslant x_2 \leqslant \varepsilon^p x_1^p  \}
  $$ 
 respectively.
 
Note that 
$$
x_1^{q} \frac{1-qr^-}{1-qs^-} \left(\frac{1- s^-}{1-r^-}\right)^{{q}}=x_2^{\frac{q}{p}} \frac{1-qr^-}{1-qs^-} \left(\frac{1-p s^-}{1-p r^-}\right)^{\frac{q}{p}}
$$
by the definition of $s^-$ and $r^-$.

Our goal is to show that $\B \leqslant B_{max}$. We will prove it using Bellman function method. The proof consists of the following parts which we will now state in the form of Lemmata.

\begin{lemma}\label{lemma1}
If the function $B_{\max}$, defined above, is concave on the domain $\Omega_\delta$, i.e. 
\begin{equation} \label{eqn:concave}
B_{\max}\left(\frac{ x^- +  x^+}{2}\right) \geqslant \frac{ B_{\max}(x^-) + B_{\max}(x^+)}{2},
\end{equation}
for any $x^+$ and $x^-$ such that there exists a weight $w \in RH_p^{\delta,d}$ with $Db^d (w) \leqslant Q$ such that $x^+ = (\av{w}{J^+}, \av{w^p}{J^+})$ and $x^- = (\av{w}{J^-}, \av{w^p}{J^-})$, then Theorem \ref{thm:Main} holds. 
\end{lemma}
\begin{lemma} \label{lemma2}
The function $B_{\max}$ is locally concave on the domain $\Omega_\varepsilon$, i.e. its Hessian matrix $d^2 B_{\max} = \left\{\frac{\partial^2 B_{\max}}{\partial x \partial y}\right\}$ is non-positive definite.
\end{lemma}
\begin{lemma}\label{lemma3}
Assume that for any three points $x$, $x^+$ and $x^- \in \Omega_\delta$ with $x = \frac{x^++x^-}{2}$, the line segment connecting $x^+$ and $x^-$ lies completely inside the larger domain $\Omega_\varepsilon$ and the
function $B_{\max}$ is locally convex on $\Omega_\varepsilon$, i.e. on $\Omega_\varepsilon$ we have that the Hessian $d^2 B_{\max}$ is non-positive definite. Then the inequality (\ref{eqn:concave}) holds for every $x^+$ and $x^- \in \Omega_\delta$.
\end{lemma}
\begin{lemma}\label{lemma4}
Let $x$, $x^+$ and $x^-$ be three points in $\Omega_\delta$ with the property that $x = \frac{x^++x^-}{2}$, then the line segment connecting $x^+$ and $x^-$ lies completely inside the larger domain $\Omega_\varepsilon$.
\end{lemma}

\subsection{Proof of Lemma \ref{lemma1}.}
First, observe that if a weight $w$ is constant on the interval $J$, $w = c$, then $\av{w^q}{J} = \av{w}{J}^q = \av{w^p}{J}^\frac{q}{p}$, therefore in this case $\B \leqslant B_{\max}$. 

Now let $w$ be a step function. Note that by concavity of $B_{\max}$ we have that since for any dyadic interval $I$ we have that $\av{w}{I} = \frac{\av{w}{I^+} + \av{w}{I^-}}{2}$ and $\av{w^p}{I} = \frac{\av{w^p}{I^+} + \av{w^p}{I^-}}{2}$
\begin{eqnarray}
|J| B_{\max} (\av{w}{J}, \av{w^p}{J}) &\geqslant& |J^-| B_{\max}(\av{w}{J^-}, \av{w^p}{J^-}) + |J^+| B_{\max}(\av{w}{J^+}, \av{w^p}{J^+})\nonumber\\ &\geqslant& |J^{--}| B_{\max}(\av{w}{J^{--}}, \av{w^p}{J^{--}}) + |J^{-+}| B_{\max}(\av{w}{J^{-+}}, \av{w^p}{J^{-+}})\nonumber\\ & & + |J^{+-}| B_{\max}(\av{w}{J^{+-}}, \av{w^p}{J^{+-}}) + |J^{++}| B_{\max}(\av{w}{J^{++}}, \av{w^p}{J^{++}}) \nonumber\\
&\geqslant& \ldots \geqslant \sum_{I\in D_n(J)} |I| B_{\max} (\av{w}{I}, \av{w^p}{I}).\nonumber
\end{eqnarray}
Now note that since $w$ is a step function, it has at most finitely many jumps. Let the number of jumps be $m$. For $n$ large enough, in the last formula we have that $w$ is constant on $2^{n}-m$ subintervals  $I\in D_n(J)$ (we will call these subintervals ``good'') and has jump discontinuities on the other $m$ subintervals (we will call them ``bad'' subintervals). On good subintervals $w$ is constant, so for such intervals we have that $|I|B_{\max} (\av{w}{I}, \av{w^p}{I}) \geqslant |I| \av{w^q}{I}$. For the bad intervals we know that $B_{\max}$ is a continuous function and the set of points $\{x = (\av{w}{I}, \av{w^p}{I}): I\in D(J) \}$ is a compact subset of $\Omega_\varepsilon$, so $B_{\max}(\av{w}{I}, \av{w^p}{I})$ for bad intervals $\{I_k\}_{k=1...m}$ are bounded by a uniform constant $M$. So the whole sum differs from $|J|\av{w^q}{J}$ by at most $mM \sum_{I ``bad''} |I|$, which tends to $0$ as $n\rightarrow \infty$. 

This implies that $\av{w}{J} \leqslant  B_{\max} (\av{w}{J}, \av{w^p}{J})$ for all step functions $w$. 

Next we extend this result to all weights $w_m$ that are bounded from above and from below $m \leqslant w \leqslant M$. We take a sequence of step-functions $w_n$ that point-wise converge to $w_m$. By the Lebesgue dominated convergence theorem Lemma \ref{lemma1} should hold for $w_m$. 

Result of \ref{Rez} extends our argument to arbitrary weight $w$, which completes the proof of the Lemma \ref{lemma1}.
\subsection{Proof of Lemma \ref{lemma2}.}

We want to show that the matrix of second derivatives of $B_{\max}$ is non-positive definite. It is not hard and it has been shown in \cite{Vasyunin:08} in a more general case. 

%

\subsection{Proof of Lemma \ref{lemma3}.}
We prove it using a nice integration trick of Nazarov, Treil and Volberg.

For the fixed points $x$, $x^+$ and $x^-$ in the domain $\Omega_\varepsilon$ such that $x = \frac{x^- + x^+}{2}$, we introduce function $b(t):= B(x_t)$, where $x_t:= \frac{1+t}{2}x^+ + \frac{1-t}{2} x^-$.  Note that defined this way, $B(x) = b(0)$ while $B(x^+) = b(1)$ and $B(x^-) = b(-1)$. Note also that 
$$
b^{\prime\prime} (t) = 
\left[ 
\begin{array}{cc}\frac{dx}{dt} & \frac{dy}{dt}
\end{array}
\right]
d^2 B_{\max}
\left[
\begin{array}{cc}
\frac{dx}{dt} \\ \frac{dy}{dt}
\end{array}
\right]
$$
So, since $-d^2 B_{\max}$ is non-negative definite, $-b^{\prime\prime} (t) \geqslant 0$ for all $-1 \leqslant t \leqslant 1$.

On the other hand, 
$$
B_{\max}(x) - \frac{B_{\max}(x^+) + B_{\max} (x^-)}{2} = b(0) - \frac{b(1) + b(-1)}{2} = -\frac{1}{2} \int_{-1}^1 (1-|t|) b^{\prime\prime} (t) dt.
$$
The second part of the above formula is a simple calculus exercise of integrating by parts twice.

Clearly, since $-b^{\prime\prime}(t)$ is non-negative, 
$$
B_{\max}(x) - \frac{B_{\max}(x^+) + B_{\max} (x^-)}{2} \geqslant 0,
$$
which completes the proof of Lemma \ref{lemma3}.
\subsection{Proof of Lemma \ref{lemma4}.}
Let $x$, $x^+$ and $x^-$ be three points in $\Omega_\delta := \{\vec{x} = (x_1, x_2) \; : \; \exists w \in RH_p^{\delta, d}\cap Db^{Q,d} \; s.t. \; x_1 = \av{w}{J}, x_2=\av{w^p}{J}\}$.
Note that the Reverse H\"older property for the weight $w$ implies that $x_1^p \leqslant x_2 \leqslant \delta^p x_1^p$ for all three points $x$, $x^+$ and $x^-$, and the fact that $w$ is almost everywhere positive implies that $x_1$, $x_2 >0$. At the same time the fact that $w$ is dyadically doubling with a doubling constant at most $Q$ implies that $x_1 \leqslant Q x_1^\pm$, $x_1^\pm \leqslant 2 x_1$ and $x^\mp_1 \leqslant (Q-1) x^\pm_1$.

Without loss of generality, we will assume that $x_1^- < x_1^+$. Then we know that $x_1 \leqslant Q x_1^-$, $x_1^+ \leqslant 2 x_1$ and $x^+_1 \leqslant (Q-1) x^-_1$.

Therefore 
$$
\Omega_\delta \subset \Omega_\delta^\prime := \{\vec{x} = (x_1, x_2) \; : \; x_1^p \leqslant x_2 \leqslant \delta^p x_1^p\},
$$
Points $x$, $x^+$ and $x^- \in \Omega_\delta^\prime$ are such that $x = \frac{x^+ +x^-}{2}$, $x^-_1 < x_1 < x_1^+$, $x_1 \leqslant Q x_1^-$, $x_1^+ \leqslant 2 x_1$ and $x^+_1 \leqslant (Q-1) x^-_1$. We need to show that the line interval connecting $x^+$ and $x^-$ lies inside the domain $\Omega_\varepsilon$.

First observe that the worst case scenario is when the central point $x$ and one of the endpoints lies on the upper boundary of $\Omega^\prime_\delta$, $x_2 = \delta^p x_1^p$, while the other endpoint lies on the lower boundary of $\Omega^\prime_\delta$, $x_2 = x_1^p$. There are two possibilities, so let us consider two cases separately.

{\bf Case 1: $x$ and $x^-$ are on the upper boundary, $x^+$ is on the lower boundary.} It means that $x = (x_1, \delta^p x_1^p)$, $x^- = (x_1^-, \delta^p (x_1^{-})^p)$ and $x^+ = (x^+_1, (x^+_1)^p)$. We need to minimize function $f(x) = x_2^{\frac{1}{p}} x_1^{-1}$ over the line that passes through $x$, $x^+$ and $x^-$. We are not going to use all conditions on our points. To simplify the problem, we will drop the condition that the point $x^+$, which is on the lower boundary. We will only be using points $x$ and $x^-$ and we will use the fact that $x_1 \leqslant Q x_1^-$.

Again, in the worst case, which may be unattainable, $x_1 = Q x_1^-$. Line through the points $x^- = (x_1^-, \delta^p(x_1^-)^p)$ and $x = (Q x^-, Q^p \delta^p (x_1^-)^p)$ has slope $\frac{\delta^p (x_1^-)^p (Q^p - 1)}{Q-1}$, therefore the equation is 
$$
y-\delta^p (x_1^-)^p - \delta^p (x_1^-)^{p-1} \frac{Q^p - 1}{Q-1} (x-x_1^-) = 0.
$$ 
So we need to solve an optimization problem
$$
\left\{
\begin{array}{cc}
f(x) = y^\frac{1}{p} x^{-1} \rightarrow \max \\
y-\delta^p (x_1^-)^p - \delta^p (x_1^-)^{p-1} \frac{Q^p - 1}{Q-1} (x-x_1^-) = 0
\end{array}
\right.
$$
The problem can be solved, for example, using method of Lagrange multipliers. If we let $H:= \frac{Q^p - 1}{Q-1}$, $f_{\max} = \left(\frac{p-1}{C-1}\right)^{\frac{p-1}{p}} \frac{C}{p} \delta$, which is exactly our choice of $\varepsilon$.

{\bf Case 2: $x$ and $x^+$ are on the upper boundary, $x^-$ is on the lower boundary.}
In this case we will also drop the condition that $x^-$ is on the lower boundary. Since coordinates of our points are positive, $x_1 = \frac{x_1^+ + x_1^-}{2} \geqslant \frac{x_1^+}{2}$, so $x_1^+ \leqslant 2 x_1$. Therefore this case is similar to the Case 1 with $Q=2$. Since $Q\geqslant 2$, this case is covered as well.

This completes the proof of Lemma \ref{lemma4} and Theorem \ref{thm:Main}.





\bibliographystyle{plain}   

\end{document}